\documentclass[12pt,a4paper]{article}
\usepackage{amsmath}
\usepackage{amssymb}
\usepackage{enumerate}
\usepackage{t1enc}
\usepackage[latin1]{inputenc}
\usepackage[english]{babel}
\usepackage{enumerate}

\usepackage{mathtools}

\DeclarePairedDelimiter\floor{\lfloor}{\rfloor}

\usepackage{amsfonts}
\usepackage{latexsym}
\usepackage{bm}
\newtheorem{theorem}{Theorem} 

\newtheorem{lemma}{Lemma}

\newtheorem{conj}{Conjecture}

\newtheorem{construction}{Construction}
\newtheorem{problem}{Problem}

\newtheorem{obs}{Observation}

\allowdisplaybreaks

\begin{document}
\title{Distance-$2$ isolation of connected graphs}

\author{Peter Borg\\[5mm]
Department of Mathematics \\
Faculty of Science \\
University of Malta\\
Malta\\
\texttt{peter.borg@um.edu.mt} 
\and 
Pawaton Kaemawichanurat\\[5mm]
Department of Mathematics \\
Faculty of Science \\
King Mongkut's University\\
of Technology Thonburi, Thailand\\
\texttt{pawaton.kae@kmutt.ac.th}
}

\date{}
\maketitle

\baselineskip = 0.53cm

\begin{abstract} 
If $\mathcal{F}$ is a set of graphs and $F$ is a copy of a graph in $\mathcal{F}$, then we call $F$ an $\mathcal{F}$-graph. Let $G$ be a connected $n$-vertex graph. For a subset $D$ of the vertex set $V(G)$, let $N^2[D]$ denote the second closed neighbourhood of $D$. Let $\iota_2(G, \mathcal{F})$ denote the size of a smallest subset $D$ of $V(G)$ such that $N^2[D]$ intersects the $\mathcal{F}$-graphs contained by $G$ (equivalently, $G - N^2[G]$ contains no $\mathcal{F}$-graph). For $k \geq 0$, let $\mathcal{E}_k$ denote the set of connected graphs that have at least $k$ edges (and at least one vertex if $k = 0$). By a well-known result of Meir and Moon, $\iota_2(G, \mathcal{E}_0) \leq n/3$ if $n \geq 3$. We show that $\iota_2(G, \mathcal{E}_1) \leq n/4$ if $n \neq 2$, $n \neq 3$ and $G$ is not a $7$-cycle. Hua and Liu recently showed that if $G$ is not a $4$-vertex path or a $4$-cycle or a $9$-cycle, then there exists a subset $D$ of $V(G)$ such that $|D| \leq n/5$ and each component of $G-N^2[D]$ has a universal vertex. We strengthen this result by showing that $\iota_2(G, \mathcal{E}_2) \leq n/5$ if $n \neq 3$, $n \neq 4$ and $G$ is not an $8$-cycle or one of five $9$-vertex graphs. It follows that $\iota_2(G, \mathcal{E}_3) \leq n/5$ if $G$ is not a $3$-cycle or a $4$-vertex graph or a $9$-cycle. This also implies the result of Hua and Liu, and that if $\mathcal{C}$ is the set of cycles, then $\iota_2(G, \mathcal{C}) \leq n/5$ unless $G$ is a $3$-cycle or a $4$-vertex graph containing a cycle. All these bounds are sharp. 
\end{abstract}

\section{Introduction}
The graph isolation problem, introduced by Caro and Hansberg in \cite{CaHa15, CaHa17}, is a natural generalization of the classical domination problem \cite{C, CH, HHS, HHS2, HL2, HL}. Similarly to the historical development of domination, the study of isolation is expanding rapidly. Many of the known isolation results are referenced in \cite{Borgrsc, Borgisdom}. The work in this paper is primarily motivated by the recent paper \cite{HuaLiu} of Hua and Liu, and the paper \cite{Borgcon} of the present first author. Inspired by the former, we start investigating a natural distance-$2$ (second closed neighbourhood) variant of the isolation problem by establishing distance-$2$ analogues of the results in \cite{Borgcon}. Two of our results strengthen the result in~\cite{HuaLiu}. 

The set of positive integers is denoted by $\mathbb{N}$. For $n \in \mathbb{N}$, $[n]$ denotes the set $\{1, \dots, n\}$ (that is, $\{i \in \mathbb{N} \colon i \leq n\}$). We take $[0]$ as the empty set $\emptyset$. Arbitrary sets are taken to be finite. For a set $X$, the set of $2$-element subsets of $X$ is denoted by ${X \choose 2}$. We may represent a $2$-element set $\{x,y\}$ by $xy$. For standard terminology in graph theory, we refer the reader to \cite{West}. Most of the terminology used here is defined in \cite{Borg1}.

Every graph $G$ is taken to be \emph{simple}, that is, $G$ is a pair $(V(G), E(G))$ such that $V(G)$ and $E(G)$ (the vertex set and the edge set of $G$) are sets that satisfy $E(G) \subseteq {V(G) \choose 2}$. We call $G$ an \emph{$n$-vertex graph} if $|V(G)| = n$. For $v \in V(G)$, $N_{G}(v)$ denotes the set of neighbours of $v$ in $G$, $N_{G}[v]$ denotes the closed neighbourhood $N_{G}(v) \cup \{ v \}$ of $v$, and $d_{G}(v)$ denotes the degree $|N_{G}(v)|$ of $v$. The maximum degree of $G$, denoted by $\Delta(G)$, is $\max\{d_G(v) \colon v \in V(G)\}$. For $X \subseteq V(G)$, $N_G[X]$ denotes the closed neighbourhood $\bigcup_{v \in X} N_G[v]$ of $X$, $N_G^2[X]$ denotes the second closed neighbourhood $N_G[N_G[X]]$ of $X$, $N_G(X)$ denotes the open neighbourhood $N_G[X] \setminus X$ of $X$, $N_G^2(X)$ denotes the second open neighbourhood $N_G^2[X] \setminus X$ of $X$, $G[X]$ denotes $(X, E(G) \cap {X \choose 2})$ (the subgraph of $G$ induced by $X$), and $G - X$ denotes $G[V(G) \setminus X]$ (the graph obtained by deleting the vertices in $X$ from $G$). We may write $N_G^2[x]$, $N_G^2(x)$ and $G-x$ instead of $N_G^2[\{x\}]$, $N_G^2(\{x\})$ and $G-\{x\}$, respectively. Where no confusion arises, the subscript $G$ may be omitted; for example, $N_G(v)$ may be abbreviated to $N(v)$.

A \emph{component of $G$} is a maximal connected subgraph of $G$. Clearly, the components of $G$ are pairwise vertex-disjoint, and their union is $G$. 

If $H$ is a copy of $G$, then we write $H \simeq G$. If $H$ is a subgraph of $G$, then we say that \emph{$G$ contains $H$}.  
  
For $n \geq 1$, the graphs $([n], {[n] \choose 2})$ and $([n], \{\{i,i+1\} \colon i \in [n-1]\})$ are denoted by $K_n$ and $P_n$, respectively. For $n \geq 3$, $C_n$ denotes the graph $([n], \{\{1,2\}, \{2,3\}, \dots, \{n-1,n\}, \{n,1\}\})$. A copy of $K_n$ is called an \emph{$n$-clique} or a \emph{complete graph}. A copy of $P_n$ is called an \emph{$n$-path} or simply a \emph{path}. A copy of $C_n$ is called an \emph{$n$-cycle} or simply a \emph{cycle}. A $3$-cycle is also called a \emph{triangle}. Note that $K_3$ is the triangle $C_3$. 

If $\mathcal{F}$ is a set of graphs and $F$ is a copy of a graph in $\mathcal{F}$, then we call $F$ an \emph{$\mathcal{F}$-graph}. A subset $D$ of $V(G)$ is called an \emph{$\mathcal{F}$-isolating set of $G$} if $N[D]$ intersects the vertex sets of the $\mathcal{F}$-graphs contained by $G$. Thus, $D$ is an $\mathcal{F}$-isolating set of $G$ if and only if $G-N[D]$ contains no $\mathcal{F}$-graph. The size of a smallest $\mathcal{F}$-isolating set of $G$ is denoted by $\iota(G, \mathcal{F})$ and called the \emph{$\mathcal{F}$-isolation number of $G$}. If $\mathcal{F} = \{F\}$, then we may replace $\mathcal{F}$ by $F$ in these terms and notation.  

Note that a subset $D$ of $V(G)$ is a $K_1$-isolating set of $G$ if and only if $N[D] = V(G)$. In the literature, such a set $D$ is called a \emph{dominating set of $G$}, and the size of a smallest dominating set of $G$ is denoted by $\gamma(G)$ and called the \emph{domination number of $G$}. Thus, $\gamma(G) = \iota(G, K_1)$. Assume that $G$ is a connected $n$-vertex graph. An early domination result of Ore \cite{Ore} is that 
\begin{equation} \gamma(G) \leq \frac{n}{2} \label{Orebound}
\end{equation}
unless $G \not\simeq K_1$ (see \cite{HHS}). Deleting the closed neighbourhood of a dominating set yields the \emph{null graph} (the graph with no vertices). Deleting the closed neighbourhood of a $K_2$-isolating set yields an \emph{empty graph} (a graph with no edges). In the literature, such a set is also called a \emph{vertex-edge dominating set} \cite{BCHH,Lewis,LHHF,Peters,TJ,Z}. Caro and Hansberg~\cite{CaHa17} proved that 
\begin{equation} \iota(G, K_2) \leq \frac{n}{3} \label{CHbound}
\end{equation}
unless $G \simeq K_2$ or $G \simeq C_5$. This solved a problem in \cite{BCHH} and was independently proved by \.{Z}yli\'{n}ski \cite{Z}. The seminal paper \cite{CaHa17} also featured several problems, some of which have been settled in \cite{Borg1,Borgrsc,BFK} (see also \cite{Borgrc2, BFK2}). As indicated above, many other isolation results followed (see \cite{Borgrsc, Borgisdom}). Domination and isolation have been particularly investigated for maximal outerplanar graphs \cite{BK, BK2, CaWa13, CaHa17, Ch75, DoHaJo16, DoHaJo17, HeKa18, LeZuZy17, Li16, MaTa96, To13, KaJi}, mostly due to connections with Chv\'{a}tal's Art Gallery Theorem \cite{Ch75}.

With a sharp upper bound established for each of $\iota(G,K_1)$ and $\iota(G,K_2)$, it was natural to consider obtaining a sharp upper bound on $\iota(G,F)$ with $F$ being any of the two connected $3$-vertex graphs $K_3$ and $P_3$, where $G$ is again a connected $n$-vertex graph. Let $\mathcal{C}$ be the set of cycles. Solving one of the problems posed by Caro and Hansberg \cite{CaHa17}, the present first author~\cite{Borg1} proved that 
\begin{equation} \iota(G,\mathcal{C}) \leq \frac{n}{4} \label{Borgcyclebound}
\end{equation} 
unless $G \simeq K_3$. Since $K_3$ is a cycle, (\ref{Borgcyclebound}) immediately implies that $\iota(G,K_3) \leq n/4$ unless $G \simeq K_3$. The significance of a $P_3$-isolating set $D$ of $G$ is that the $E(G-N[D])$ is a matching, that is, no two edges of $G-N[D]$ have a common vertex (so $\Delta(G-N[D]) \leq 1$). Zhang and Wu \cite{ZW} proved that  
\begin{equation} \iota(G,P_3) \leq \frac{2n}{7} \label{ZWbound}
\end{equation} 
unless $G$ is a $\{P_3, C_3, C_6\}$-graph. We mention in passing that the present first author \cite{Borgcon} established the stronger inequality $\iota(G,P_3) \leq (4n-r)/14$ with $r$ being the number of leaves of $G$, and $G$ not being a copy of one of six particular graphs. 
The $P_3$-isolation problem attracted further attention and was consequently also treated in \cite{BBS, BBS2, CLWX, HZJ}.

We now introduce a variant of $\mathcal{F}$-isolation that generalizes distance-$2$ domination. If $G$ is a graph and $D \subseteq V(G)$ such that $N^2[D]$ intersects the vertex sets of the $\mathcal{F}$-graphs contained by $G$, then we call $D$ a \emph{distance-$2$ $\mathcal{F}$-isolating set of $G$} or, in brief, an \emph{$(\mathcal{F}, 2)$-isolating set of $G$}. Thus, $D$ is an $(\mathcal{F}, 2)$-isolating set of $G$ if and only if $G-N^2[D]$ contains no $\mathcal{F}$-graph. The size of a smallest $(\mathcal{F}, 2)$-isolating set of $G$ is denoted by $\iota_2(G, \mathcal{F})$ and called the \emph{distance-2 $\mathcal{F}$-isolation number of $G$} or, in brief, the \emph{$(\mathcal{F}, 2$)-isolation number of $G$}. If $\mathcal{F} = \{F\}$, then we may replace $\mathcal{F}$ by $F$ in these terms and notation. In the literature, the $(K_1,2)$-isolation number of $G$ is called the \emph{distance-$2$ domination} (or \emph{distance $2$-domination}) \emph{number of $G$}.

The following is a sharp bound that follows from a well-known result of Meir and Moon \cite{MM} and was also proved in \cite{SSE}.

\begin{theorem}[\cite{MM, SSE}] \label{result k=0}
If $G$ is a connected $n$-vertex graph with $n \geq 3$, then
\[\iota_2(G,K_1) \leq \left \lfloor \frac{n}{3} \right \rfloor .\] 
\end{theorem}

Let $\mathcal{E}_0$ denote the set of connected graphs with non-empty vertex sets. For $k \geq 1$, let $\mathcal{E}_k$ denote the set of connected graphs that have at least $k$ edges. For $j \geq k$, $\iota_2(G,\mathcal{E}_j) \leq \iota_2(G, \mathcal{E}_k)$ as every $(\mathcal{E}_k, 2)$-isolating set of $G$ is an $(\mathcal{E}_{j}, 2)$-isolating set of $G$. It is easy to see that, similarly, \\
- a $(K_1, 2)$-isolating set is an $(\mathcal{E}_0, 2)$-isolating set, and vice-versa, \\
- a $(K_2, 2)$-isolating set is an $(\mathcal{E}_1, 2)$-isolating set, and vice-versa, and \\
- a $(P_3, 2)$-isolating set is an $(\mathcal{E}_2, 2)$-isolating set, and vice-versa. \\
Therefore,
\begin{equation} \iota_2(G,\mathcal{E}_0) = \iota_2(G, K_1), \quad \quad  \iota_2(G,\mathcal{E}_1) = \iota_2(G, K_2), \quad \quad \iota_2(G,\mathcal{E}_2) = \iota_2(G, P_3). \label{isoleqns}
\end{equation}
In this paper, we introduce the problem of determining a sharp upper bound on $\iota_2(G, \mathcal{E}_k)$ similar to the bounds above. Theorem~\ref{result k=0} is the solution for $k = 0$. We solve the problem for $k = 1, 2, 3$. The following construction will be used to show that our bounds are sharp.

\begin{construction} \label{const} \emph{For any $n, k \in \mathbb{N}$, we construct two connected $n$-vertex graphs $B_{n,k}$ and $B_{n,k}'$ as follows.  If $n \leq k$, then let $B_{n,k} = B_{n,k}' = P_n$. Suppose $n \geq k+1$. Let $a = \left \lfloor \frac{n}{k+1} \right \rfloor$ and $b = n - ka$ (so $a \leq b \leq a+k$). Let $F_1, \dots, F_a$ be copies of $P_k$ such that $P_b, F_1, \dots, F_a$ are pairwise vertex-disjoint. For each $i \in [a]$, let $v_{i,1}, \dots, v_{i,k}$ be the vertices of $F_i$ corresponding to the vertices $1, \dots, k$ of $P_k$, respectively. Let $B_{n,k}$ be the graph with $V(B_{n,k}) = [b] \cup \bigcup_{i=1}^a V(F_i)$ and $E(B_{n,k}) = E(P_a) \cup \{\{a, j\} \colon j \in [b] \setminus [a]\} \cup \bigcup_{i=1}^a (E(F_i) \cup \{\{i,v_{i,1}\}\})$. If $k \leq 2$, then let $B_{n,k}' = B_{n,k}$. If $k \geq 3$, then let $B_{n,k}' = (V(B_{n,k}), E(B_{n,k}) \cup \{v_{i,k-2}v_{i,k} \colon i \in [a]\})$.}
\end{construction}

It is easy to see that if $G = B_{n,2}$ with $n \geq 3$, then $G$ attains the bound in Theorem~\ref{result k=0} (see, for example, the proof of Theorem~\ref{result k=1}).

We may abbreviate $\iota_2(G,K_2)$ to $\iota_2(G)$. If $G$ is a connected $n$-vertex graph, then $\iota_2(G) = 1 = \frac{n+2}{4}$ if $n = 2$, and $\iota_2(G) = \frac{n+1}{4}$ if $n = 3$ or $G \simeq C_7$.  In Section~\ref{Proofsection}, we prove the following result, which, by (\ref{isoleqns}), settles our problem for $k = 1$.

\begin{theorem} \label{result k=1}
If $G$ is a connected $n$-vertex graph, $n \notin \{2, 3\}$ and $G$ is not a $7$-cycle, then
\[\iota_2(G) \leq \left \lfloor \frac{n}{4} \right \rfloor .\] 
Moreover, equality holds if $G = B_{n,3}$.
\end{theorem}

For each $i \in \{3, 4\}$, let $\mathcal{G}_i$ be the set of connected graphs whose vertex set is $[i]$. Let $\mathcal{G}_8 = \{C_8\}$. Let $G_9^1 = ([9], E(C_8) \cup \{\{1,9\}\})$, $G_9^2 = ([9], E(G_9^1) \cup \{\{2,9\}\})$, $G_9^3 = ([9], E(G_9^1) \cup \{\{3,9\}\})$ and $G_9^4 = ([9], E(G_9^2) \cup \{\{3,9\}\})$. Let $\mathcal{G}_9 = \{C_9, G_9^1, G_9^2, G_9^3, G_9^4\}$. Let $\mathcal{G} = \mathcal{G}_3 \cup \mathcal{G}_4 \cup \mathcal{G}_8 \cup \mathcal{G}_9$. It is easy to see that if $G$ is a $\mathcal{G}_i$-graph for some $i \in \{3,4,8,9\}$, then $\iota_2(G,P_3) = \frac{i+2}{5}$ if $i \in \{3,8\}$, and $\iota_2(G,P_3) = \frac{i+1}{5}$ if $i \in \{4,9\}$. In Section~\ref{Proofsection2}, we prove the following result, which, by (\ref{isoleqns}), settles our problem for $k = 2$.
\begin{theorem} \label{result k=2}
If $G$ is a connected $n$-vertex graph that is not a $\mathcal{G}$-graph, then
\[\iota_2(G, P_3) \leq \left \lfloor \frac{n}{5} \right \rfloor .\] 
Moreover, equality holds if $G = B_{n,4}$.
\end{theorem}
Theorem~\ref{result k=2} is intriguing in more than one way. Firstly, the proof uses almost all the ideas produced to date, such as multiple vertex deletion steps and divisibility considerations, for establishing such isolation bounds, indicating that distance variants of isolation may necessitate novel approaches. Secondly, it yields our next result, which tells us that the best possible bound on $\iota_2(G, \mathcal{E}_3)$ is the same as the one on $\iota_2(G, \mathcal{E}_2)$ (since $\iota_2(G, \mathcal{E}_2) = \iota_2(G, P_3)$ by (\ref{isoleqns})).
\begin{theorem} \label{result k=3}
If $G$ is a connected $n$-vertex graph, $n \neq 4$ and $G$ is not a $\{C_3, C_9\}$-graph, then
\[\iota_2(G, \mathcal{E}_3) \leq \left \lfloor \frac{n}{5} \right \rfloor .\] 
Moreover, equality holds if $G = B_{n,4}'$.
\end{theorem}
\textbf{Proof.} As pointed out above, $\iota_2(G,\mathcal{E}_3) \leq \iota_2(G,\mathcal{E}_2)$, so the bound follows by Theorem~\ref{result k=2}. Suppose that $G$ is the graph $B_{n,k}'$ in Construction~\ref{const} with $k = 4$. Note that $[a]$ is an $\mathcal{E}_3$-isolating set of $G$, so $\iota_2(G,\mathcal{E}_3) \leq a$. For each $i \in [a]$, $F_i - v_{i,1} \simeq K_3$ and $N^2(V(F_i - v_{i,1})) = \{v_{i,1}, i\}$. Thus, if $D$ is an $(\mathcal{E}_3,2)$-isolating set of $G$, then $D \cap (V(F_i) \cup \{i\}) \neq \emptyset$ for each $i \in [a]$. Therefore, $\iota_2(G,\mathcal{E}_3) = a$.~\hfill{$\Box$}
\\

Theorem~\ref{result k=3} in turn yields the following analogue of (\ref{Borgcyclebound}).
\begin{theorem} \label{result cycles}
If $G$ is a connected $n$-vertex graph that is not a $3$-cycle or a $4$-vertex graph containing a cycle, then
\[\iota_2(G, \mathcal{C}) \leq \left \lfloor \frac{n}{5} \right \rfloor .\] 
Moreover, equality holds if $G = B_{n,4}'$.
\end{theorem}
Since $\mathcal{C} \subset \mathcal{E}_3$, $\iota_2(G,\mathcal{C}) \leq \iota_2(G,\mathcal{E}_3)$, so the bound follows by Theorem~\ref{result k=3}. The second part of the theorem is obtained as in the proof of Theorem~\ref{result k=3}. Hua and Liu \cite{HuaLiu} recently showed that if $G$ is not a $\{P_4, C_4, C_9\}$-graph, then $V(G)$ has a subset $D$ such that $|D| \leq n/5$ and each component of $G-N^2[D]$ has a universal vertex (a vertex adjacent to all the other vertices). This follows immediately from each of Theorems~\ref{result k=2} and \ref{result k=3}.

\section{Proof of Theorem~\ref{result k=1}} \label{Proofsection}

In this section, we prove Theorem~\ref{result k=1}. We start with analogues of \cite[Lemmas~2 and 3]{Borg1} that will be used repeatedly. 

\begin{lemma} \label{lemma}
If $G$ is a graph, $\mathcal{F}$ is a set of graphs, $X \subseteq V(G)$ and $Y \subseteq N^2[X]$, then \[\iota_2(G, \mathcal{F}) \leq |X| + \iota_2(G-Y, \mathcal{F}).\] 
\end{lemma}
\textbf{Proof.} Let $D$ be an $(\mathcal{F}, 2)$-isolating set of $G-Y$ of size $\iota_{2}(G-Y, \mathcal{F})$. Clearly, $\emptyset \neq V(F) \cap Y \subseteq V(F) \cap N^2[X]$ for each $\mathcal{F}$-graph $F$ that is a subgraph of $G$ and not a subgraph of $G-Y$. Thus, $X \cup D$ is an $(\mathcal{F}, 2)$-isolating set of $G$. The result follows.~\hfill{$\Box$}
\\

For a graph $G$, let ${\rm C}(G)$ denote the set of components of $G$.

\begin{lemma} \label{lemmacomp}
If $G$ is a graph and $\mathcal{F}$ is a set of connected graphs, then \[\iota_2(G, \mathcal{F}) = \sum_{H \in {\rm C}(G)} \iota_2(H, \mathcal{F}).\] 
\end{lemma}
\textbf{Proof.} For each $H \in {\rm C}(G)$, let $D_H$ be a smallest $(\mathcal{F}, 2)$-isolating set of $H$. Consider any $\mathcal{F}$-graph $F$ contained by $G$. Since $F$ is connected, $F$ is contained by some $I \in {\rm C}(G)$, so $N_I^2[D_I] \cap V(F) \neq \emptyset$. Therefore, $\bigcup_{H \in {\rm C}(G)} D_H$ is an $(\mathcal{F}, 2)$-isolating set of $G$, and hence $\iota_2(G,\mathcal{F}) \leq \sum_{H \in {\rm C}(G)} |D_H| = \sum_{H \in {\rm C}(G)} \iota_2(H,\mathcal{F})$. Now let $D$ be a smallest $(\mathcal{F}, 2)$-isolating set of $G$. For each $H \in {\rm C}(G)$, $D \cap V(H)$ is an $(\mathcal{F}, 2)$-isolating set of $H$. We have $\sum_{H \in {\rm C}(G)} \iota_2(H,\mathcal{F}) \leq \sum_{H \in {\rm C}(G)} |D \cap V(H)| = |D| = \iota_2(G,\mathcal{F}) \leq \sum_{H \in {\rm C}(G)} \iota_2(H,\mathcal{F})$, so equality holds throughout.~\hfill{$\Box$} 
\\

If $n \geq 4$, then $\{4i \colon 1 \leq i \leq \lfloor n/4 \rfloor\}$ is a $(K_2, 2)$-isolating set of $P_n$ of size $\lfloor n/4 \rfloor$, and $\{6i-5 \colon 1 \leq i \leq \lfloor (n+5)/6 \rfloor\}$ is a $(K_2, 2)$-isolating set of $C_n$ of size $\lfloor (n+5)/6 \rfloor$. If $n \geq 10$, then $(n+5)/6 \leq n/4$. If $4 \leq n \leq 9$ and $n \neq 7$, then $\lfloor (n+5)/6 \rfloor \leq \lfloor n/4 \rfloor$. Consequently, 
\begin{equation} \mbox{if $n \geq 4$, $G \in \{P_n, C_n\}$ and $G \neq C_7$, then $\iota_2(G, K_2) \leq \lfloor n/4 \rfloor$.}\medskip \label{pcbound2}
\end{equation}

\noindent
\textbf{Proof of Theorem~\ref{result k=1}.} Suppose that $G$ is the graph $B_{n,k}$ in Construction~\ref{const} with $k = 3$. Note that $[a]$ is a $(K_2, 2)$-isolating set of $G$, so $\iota_2(G) \leq a$. For each $i \in [a]$, $F_i - v_{i,1} \simeq K_2$ and $N^2(V(F_i - v_{i,1})) = \{v_{i,1}, i\}$. Thus, if $D$ is a $(K_2,2)$-isolating set of $G$, then $D \cap (V(F_i) \cup \{i\}) \neq \emptyset$ for each $i \in [a]$. Therefore, $\iota_2(G) = a$.

We now prove that the bound in the theorem holds, using induction on $n$. Since $\iota_2(G)$ is an integer, it suffices to prove that $\iota_2(G) \leq n/4$. If $n = 1$, then trivially $\iota_2(G) = 0 = \lfloor n/4 \rfloor$. Suppose $n > 1$. Since $n \notin \{2,3\}$, $n \geq 4$. Since $G$ is connected, if $n = 4$, then for each $v \in V(G)$, $\{v\}$ is a $(K_2, 2)$-isolating set of $G$, so $\iota_2(G) = 1 = \floor{n/4}$. Suppose $n \geq 5$. If $\Delta(G) \leq 2$, then $G$ is a path or a cycle, and hence the result is given by (\ref{pcbound2}). Suppose $\Delta(G) \geq 3$. Let $v \in V(G)$ such that $d(v) = \Delta(G)$. Then, $|N[v]| \geq 4$. If $V(G) = N[v]$, then $\{v\}$ is a $(K_2,2)$-isolating set of $G$, so $\iota_2(G) = 1 \leq n/4$. Suppose $V(G) \neq N[v]$. Let $G' = G-N[v]$. Then, $V(G') \neq \emptyset$. Let $\mathcal{H} = {\rm C}(G')$. Let 
\begin{equation} \mathcal{H}' = \{H \in \mathcal{H} \colon H \mbox{ is a $\{K_2, P_3, C_3, C_7\}$-graph}\}. \nonumber 
\end{equation} 
Let $\mathcal{H}^* = \mathcal{H} \setminus \mathcal{H}'$. For each $H \in \mathcal{H}^*$, let $D_H$ be a smallest $(K_2, 2)$-isolating set of $H$ of size $\iota_2(H)$. By the induction hypothesis, $|D_H| \leq n/4$ for each $H \in \mathcal{H}^*$. If $\mathcal{H}' = \emptyset$, then by Lemma~\ref{lemma} (with $X = \{v\}$ and $Y = N[v]$) and Lemma~\ref{lemmacomp},
\begin{equation} \iota_2(G) \leq 1 + \iota_2(G') = 1 + \sum_{H \in \mathcal{H}^*} \iota_2(H) \leq \frac{|N[v]|}{4} + \sum_{H \in \mathcal{H}^*} \frac{|V(H)|}{4} = \frac{n}{4}. \nonumber 
\end{equation}
Suppose $\mathcal{H}' \neq \emptyset$. 

For any $H \in \mathcal{H}$ and any $x \in N(v)$ such that $xy_{x,H} \in E(G)$ for some $y_{x,H} \in V(H)$, we say that $H$ is \emph{linked to $x$} and that $x$ is \emph{linked to $H$}. Since $G$ is connected, each member of $\mathcal{H}$ is linked to at least one member of $N(v)$. Thus, for each $H \in \mathcal{H}$, $x_Hy_H \in E(G)$ for some $x_H \in N(v)$ and $y_H \in V(H)$.
\\
\\
\textbf{Case 1:} \emph{$I \simeq C_7$ for some $I \in \mathcal{H}'$}. Then, $I$ has a vertex $y_I'$ such that $y_I \notin N_I^2[y_I']$. Let $Y = N_I^2[y_I']$. Let $G_I = G - Y$. Since $I - Y \simeq K_2$ and $y_I \in V(I-Y)$, $G_I$ is connected. We have $N_{G_I}[v] = N[v]$ and $|V(G_I)| \geq |N[v]| + |V(I-Y)| \geq 6$. Since $d_{G_I}(v) \geq 3$, $G_I \not\simeq C_7$. By the induction hypothesis, $\iota_2(G_I) \leq |V(G_I)|/4 = (n-5)/4$. By Lemma~\ref{lemma}, $\iota_2(G) \leq 1 + \iota_2(G_I) < n/4$.
\\
\\
\textbf{Case 2:} \emph{$H \not\simeq C_7$ for each $H \in \mathcal{H}'$}. Then, each member of $\mathcal{H}'$ is a $\{K_2, P_3, C_3\}$-graph. For each $i \in \{2, 3\}$, let $\mathcal{H}_i = \{H \in \mathcal{H}' \colon |V(H)| = i\}$. Thus, $\mathcal{H}' = \mathcal{H}_2 \cup \mathcal{H}_3$. Let $Y = \{y_H \colon H \in \mathcal{H}_2\}$ and $Y' = N[v] \cup Y$. We have $Y' \subseteq N^2[v]$ and ${\rm C}(G-Y') = \mathcal{H}^* \cup \{H - y_H \colon H \in \mathcal{H}_2\} \cup \mathcal{H}_3$.
\\
\\
\textbf{Case 2.1:} $\mathcal{H}_3 = \emptyset$. Then, $\mathcal{H}' = \mathcal{H}_2$. By Lemmas~\ref{lemma} and \ref{lemmacomp}, 
\begin{align} \iota_2(G) &\leq 1 + \iota_2(G-Y') = 1 + \sum_{H \in \mathcal{H}^*} \iota_2(H) + \sum_{H \in \mathcal{H}_2} \iota_2(H - y_H) \nonumber \\
&\leq \frac{|N[v]|}{4} + \sum_{H \in \mathcal{H}^*} \frac{|V(H)|}{4} + \sum_{H \in \mathcal{H}_2} \iota_2(K_1) < \frac{n}{4} \nonumber
\end{align}
as $\mathcal{H}' \neq \emptyset$ and $\iota_2(K_1) = 0$.
\\
\\
\textbf{Case 2.2:} $\mathcal{H}_3 \neq \emptyset$. Let $X = \{x_H \colon H \in \mathcal{H}_3\}$. Then, $|X| \leq |\mathcal{H}_3|$. We have $4|X| \leq |X| + 3|\mathcal{H}_3| = |X| + \sum_{H \in \mathcal{H}_3} |V(H)|$, so
\begin{equation} |X| \leq \frac{|X|}{4} + \sum_{H \in \mathcal{H}_3} \frac{|V(H)|}{4}. \label{4|X|}
\end{equation}

Suppose $x_H \in X$ for each $H \in \mathcal{H}_2$. Let $X' = N[v] \cup \bigcup_{H \in \mathcal{H}'} (N^2[x_H] \cap V(H))$. Then, $X' \subseteq N^2[X]$, $\mathcal{H}^* \subset {\rm C}(G-X')$ and $H \simeq K_1$ for each $H \in {\rm C}(G-X') \setminus \mathcal{H}^*$. By Lemma~\ref{lemma}, Lemma~\ref{lemmacomp} and (\ref{4|X|}), 
$$\iota_2(G) \leq |X| + \iota_2(G-X') = |X| + \sum_{H \in \mathcal{H}^*} \iota_2(H) \leq \frac{|X|}{4} + \sum_{H \in \mathcal{H}_3} \frac{|V(H)|}{4} + \sum_{H \in \mathcal{H}^*} \frac{|V(H)|}{4} < \frac{n}{4}.$$

Now suppose $x_I \notin X$ for some $I \in \mathcal{H}_2$. Let $X^+ = \{v\} \cup X$ and $X' = \{v\} \cup Y \cup \bigcup_{H \in \mathcal{H}_3} (N^2[x_H] \cap V(H))$. Then, $X' \subseteq N^2[X^+]$, $\mathcal{H}^* \subset {\rm C}(G-X')$ and $H \simeq K_1$ for each $H \in {\rm C}(G-X') \setminus \mathcal{H}^*$. By Lemma~\ref{lemma}, Lemma~\ref{lemmacomp} and (\ref{4|X|}),  
\begin{align} \iota_2(G) &\leq |X^+| + \iota_2(G-X') = 1 + |X| + \sum_{H \in \mathcal{H}^*} \iota_2(H) \nonumber \\
&\leq \frac{|\{v,x_I\} \cup V(I)|}{4} + \frac{|X|}{4} + \sum_{H \in \mathcal{H}_3} \frac{|V(H)|}{4} + \sum_{H \in \mathcal{H}^*} \frac{|V(H)|}{4} \leq \frac{n}{4}, \nonumber
\end{align}
as required.\hfill{$\Box$}

\section{Proof of Theorem~\ref{result k=2}} \label{Proofsection2}

We now start working towards the proof of Theorem~\ref{result k=2}.

\begin{lemma}\label{5vertex} If $G$ is a connected $5$-vertex graph, then $N^2[v] = V(G)$ for some $v \in V(G)$. 
\end{lemma}
\textbf{Proof.} If $\Delta(G) \leq 2$, then $G$ is a $5$-path or a $5$-cycle, so the result is straightforward. Suppose $\Delta(G) \geq 3$. Let $v \in V(G)$ with $d(v) = \Delta(G)$. If $d(v) = 4$, then $N[v] = V(G)$. Suppose $d(v) = 3$. Let $w \in V(G) \setminus N[v]$. Since $|V(G)| = 5 = |N[v]| + 1$, $V(G) \setminus N[v] = \{w\}$. Since $G$ is connected, $w \in N(x)$ for some $x \in N(v)$, so $N^2[v] = V(G)$.\hfill{$\Box$}
\\

We will make use of the following straightforward observation.

\begin{obs}\label{8cycleobs} If $i \in [4]$, $G \simeq G_9^i$ and $v \in V(G)$ with $d(v) = 3$, then $G$ contains an $8$-cycle $F = (\{v, x_1, \dots, x_7\}, \{vx_1, x_1x_2, \dots, x_6x_7, x_7v\})$ and has a vertex $x_8$ such that $\{x_8\} = V(G) \setminus V(F)$ and $vx_8 \in E(G) \setminus E(F) \subseteq \{vx_8, x_1x_8, x_2x_8\}$.
\end{obs}

If $n \geq 5$, then $\{5i \colon 1 \leq i \leq \lfloor n/5 \rfloor\}$ is a $(P_3, 2)$-isolating set of $P_n$ of size $\lfloor n/5 \rfloor$, and $\{7i-6 \colon 1 \leq i \leq \lfloor (n+6)/7 \rfloor\}$ is a $(P_3, 2)$-isolating set of $C_n$ of size $\lfloor (n+6)/7 \rfloor$. If $n \geq 15$, then $(n+6)/7 \leq n/5$. If $5 \leq n \leq 14$ and $n \notin \{8, 9\}$, then $\lfloor (n+6)/7 \rfloor \leq \lfloor n/5 \rfloor$. Consequently, 
\begin{equation} \mbox{if $n \geq 5$, $G \in \{P_n, C_n\}$ and $G \notin \{C_8, C_9\}$, then $\iota_2(G, P_3) \leq \lfloor n/5 \rfloor$.}\medskip \label{pcbound}
\end{equation}

\noindent \textbf{Proof of Theorem~\ref{result k=2}.} Suppose that $G$ is the graph $B_{n,k}$ in Construction~\ref{const} with $k = 4$. Note that $[a]$ is a $(P_3,2)$-isolating set of $G$, so $\iota_2(G, P_3) \leq a$. For each $i \in [a]$, $F_i - v_{i,1} \simeq P_3$ and $N^2(V(F_i - v_{i,1})) = \{v_{i,1}, i\}$. Thus, if $D$ is a $(P_3,2)$-isolating set of $G$, then $D \cap (V(F_i) \cup \{i\}) \neq \emptyset$ for each $i \in [a]$. Therefore, $\iota_2(G, P_3) = a$.

We now prove that the bound in the theorem holds, using induction on $n$. Since $\iota_2(G,P_3)$ is an integer, it suffices to prove that $\iota_2(G) \leq n/5$. If $n \leq 2$, then trivially $\iota_2(G,P_3) = 0 = \lfloor n/5 \rfloor$. Suppose $n > 2$. Since $n \notin \{3,4\}$, $n \geq 5$. If $\Delta(G) \leq 2$, then $G$ is a path or a cycle, and hence the result is given by (\ref{pcbound}). Suppose $\Delta(G) \geq 3$. Let $v \in V(G)$ such that $d(v) = \Delta(G)$. Let $G' = G - N[v]$. Then, $|N[v]| \geq 4$. If $V(G) = N[v]$, then $\{v\}$ is a $(P_3,2)$-isolating set of $G$, so $\iota_2(G,P_3) = 1 \leq n/5$. Suppose $V(G) \neq N[v]$. Let $G' = G-N[v]$. Then, $V(G') \neq \emptyset$. Let $\mathcal{H} = {\rm C}(G')$. As in the proof of Theorem~\ref{result k=1}, for each $H \in \mathcal{H}$, $H$ is linked to at least one member of $N(v)$, so $x_Hy_H \in E(G)$ for some $x_H \in N(v)$ and $y_H \in V(H)$.
 
For each $i \in \{3,4,8,9\}$, let $\mathcal{H}_i = \{H \in \mathcal{H} \colon H \text{ is a $\mathcal{G}_i$-graph}\}$. Let 
$$\mathcal{H}' = \mathcal{H}_3 \cup \mathcal{H}_4 \cup \mathcal{H}_8 \cup \mathcal{H}_9.$$ 
Let $\mathcal{H}^* = \mathcal{H} \setminus \mathcal{H}'$. For each $H \in \mathcal{H}^*$, let $D_H$ be a smallest $(P_3, 2)$-isolating set of $H$ of size $\iota_2(H, P_3)$. By the induction hypothesis, $|D_H| \leq n/5$ for each $H \in \mathcal{H}^*$. 

Suppose $\mathcal{H}_8 \cup \mathcal{H}_9 \neq \emptyset$. Let $I \in \mathcal{H}_8 \cup \mathcal{H}_9$. We may assume that $I \in \mathcal{G}_8 \cup \mathcal{G}_9$. Let 
$$y_I' = \begin{cases}
    y_I + 4 \quad &\mbox{if $y_I \in [4]$,} \\
    y_I - 4 \quad &\mbox{if $y_I \in [9] \setminus [4]$.}
\end{cases}$$
Let $Y = N_I^2[y_I']$ and $G_I = G - Y$. 
Then, $y_I \in V(I - Y)$, $|N_{I-Y}[y_I]| \geq 3$, $G_I$ has a component $G_v^*$ such that $N[v] \cup N_{I-Y}[y_I] \subseteq V(G_v^*)$, and ${\rm C}(G_I) \subseteq \{G_v^*, (\{9\}, \emptyset)\}$. We have $|V(G_v^*)| \geq 7$. By Lemmas~\ref{lemma} and \ref{lemmacomp}, $\iota_2(G,P_3) \leq 1 + \iota_2(G_I,P_3) = 1 + \iota_2(G_v^*,P_3)$. If $|V(G_v^*)| \notin \{8, 9\}$, then by the induction hypothesis, $\iota_2(G_v^*,P_3) \leq |V(G_v^*)|/5 \leq (n-5)/5$, so $\iota_2(G,P_3) \leq n/5$. If $|V(G_v^*)| \in \{8, 9\}$, then since $N[v] \cup N_{I-Y}[y_I]$ is a subset of $N_{G_v^*}^2[x_I]$ of size at least $7$, $\{x_I\}$ is a $(P_3,2)$-isolating set of $G_v^*$, so $\iota_2(G,P_3) < n/5$.

Now suppose $\mathcal{H}_8 \cup \mathcal{H}_9 = \emptyset$. Then, $\mathcal{H}' = \mathcal{H}_3 \cup \mathcal{H}_4$. Let $X = \{x_H \colon H \in \mathcal{H}_4\}$. Then, $|X| \leq |\mathcal{H}_4|$. We have $5|X| \leq |X| + 4|\mathcal{H}_4| = |X| + \sum_{H \in \mathcal{H}_4}|V(H)|$, so 
\begin{equation} |X| \leq \frac{|X|}{5} + \sum_{H \in \mathcal{H}_4} \frac{|V(H)|}{5}. \label{5|X|}
\end{equation}
\medskip

\noindent \textbf{Case 1:} $\mathcal{H}_3 \neq \emptyset$. 

Suppose $X = N(v)$. Let $X' = N[v] \cup \bigcup_{H \in \mathcal{H}'} (N^2[x_H] \cap V(H))$. Then, $X' \subseteq N^2[X]$, $\mathcal{H}^* \subset {\rm C}(G-X')$ and $|V(H)| \leq 2$ for each $H \in {\rm C}(G-X') \setminus \mathcal{H}^*$. By Lemma~\ref{lemma}, Lemma~\ref{lemmacomp} and (\ref{5|X|}), 
\begin{align} \iota_2(G,P_3) &\leq |X| + \iota_2(G-X',P_3) = |X| + \sum_{H \in \mathcal{H}^*} \iota_2(H,P_3) \nonumber \\
&\leq \frac{|X|}{5} + \sum_{H \in \mathcal{H}_4} \frac{|V(H)|}{5} + \sum_{H \in \mathcal{H}^*} \frac{|V(H)|}{5} < \frac{n}{5}. \nonumber
\end{align}

Now suppose $X \neq N(v)$. Let $Z = \{v\} \cup (N(v) \setminus X) \cup \bigcup_{H \in \mathcal{H}_3} V(H)$. Since $\mathcal{H}_3 \neq\emptyset$ and $X \subsetneq N(v)$, $|Z| \geq 5$. Let $X^+ = \{v\} \cup X$ and 
$$X' = N[v] \cup \{y_H \colon H \in \mathcal{H}_3\} \cup \bigcup_{H \in \mathcal{H}_4} (N^2[x_H] \cap V(H)).$$
Then, $X' \subseteq N^2[X^+]$, $\mathcal{H}^* \subset {\rm C}(G-X')$ and $|V(H)| \leq 2$ for each $H \in {\rm C}(G-X') \setminus \mathcal{H}^*$. By Lemma~\ref{lemma}, Lemma~\ref{lemmacomp} and (\ref{5|X|}),  
\begin{align} \iota_2(G,P_3) &\leq |X^+| + \iota_2(G-X',P_3) = 1 + |X| + \sum_{H \in \mathcal{H}^*} \iota_2(H,P_3) \nonumber \\
&\leq \frac{|Z|}{5} + \frac{|X|}{5} + \sum_{H \in \mathcal{H}_4} \frac{|V(H)|}{5} + \sum_{H \in \mathcal{H}^*} \frac{|V(H)|}{5} \leq \frac{n}{5}. \nonumber
\end{align}
\medskip

\noindent \textbf{Case 2:} $\mathcal{H}_3 = \emptyset$. Then, $\mathcal{H}' = \mathcal{H}_4$ and $X = \{x_H \colon H \in \mathcal{H}'\}$. If $\mathcal{H}' \neq \emptyset$, then we obtain $\iota_2(G, P_3) < n/5$ as in the case $X = N(v)$ of Case~1. Suppose $\mathcal{H}' = \emptyset$. Then, 
$$\mathcal{H} = \mathcal{H}^*.$$ 
If $|N[v]| \geq 5$, then by Lemmas~\ref{lemma} and \ref{lemmacomp},
$$\iota_2(G, P_3) \leq 1 + \iota_2(G', P_3) \leq \frac{|N[v]|}{5} + \sum_{H \in \mathcal{H}^*}\frac{|V(H)|}{5} \leq  \frac{n}{5}.$$
Suppose $|N[v]| \leq 4$. We have $3 \leq \Delta(G) = d(v) \leq 3$, so $\Delta(G) = 3$. Let $I \in \mathcal{H}^*$, $I' = I - y_I$ and $\mathcal{I} = {\rm C}(I')$. Let $Y_I = N[v] \cup \{y_I\}$, $G_I = G - Y_I$ and $\mathcal{H}_I = (\mathcal{H}^* \setminus \{I\}) \cup \mathcal{I}$. We have $Y_I \subseteq N^2[x_I]$. For each $i \in \{3, 4, 8, 9\}$, let $\mathcal{I}_i = \{H \in \mathcal{I} \colon H \text{ is a  $\mathcal{G}_i$-graph}\}$. Let $\mathcal{I}' = \mathcal{I}_3 \cup \mathcal{I}_4 \cup \mathcal{I}_8 \cup \mathcal{I}_9$ and $\mathcal{I}^* = \mathcal{I} \setminus \mathcal{I}'$. For each $H \in \mathcal{I}^*$, let $D_H$ be a smallest $(P_3, 2)$-isolating set of $H$ of size $\iota_2(H, P_3)$. By the induction hypothesis, $|D_H| \leq n/5$ for each $H \in \mathcal{I}^*$. If $\mathcal{I} = \mathcal{I}^*$, then by Lemmas~\ref{lemma} and \ref{lemmacomp}, 
$$\iota_2(G, P_3) \leq 1 + \iota_2(G_I, P_3) = \frac{|Y_I|}{5} + \sum_{H \in \mathcal{H}_I} \iota_2(G, P_3) \leq \frac{|Y_I|}{5} + \sum_{H \in \mathcal{H}_I}\frac{|V(H)|}{5} = \frac{n}{5}.$$
Suppose $\mathcal{I} \neq \mathcal{I}^*$. Thus, 
$$\mathcal{I}' \neq \emptyset.$$ 

Since $I$ is connected, for each $H \in \mathcal{I}$, $y_Iz_H \in E(G)$ for some $z_H \in V(H)$. Since $\Delta(G) = 3$ and $x_I \in N(y_I)$, it follows that $|\{z_H \colon H \in \mathcal{I}\}| \leq 2$, so
$$|\mathcal{I}| \leq 2.$$ 

Suppose $\mathcal{I}_8 \cup \mathcal{I}_9 \neq \emptyset$. Let $J \in \mathcal{I}_8 \cup \mathcal{I}_9$. We may assume that $J \in \mathcal{G}_8 \cup \mathcal{G}_9$. Let 
$$z_J' = \begin{cases}
    z_J + 4 \quad &\mbox{if $z_J \in [4]$,} \\
    z_J - 4 \quad &\mbox{if $z_J \in [9] \setminus [4]$.}
\end{cases}$$
Let $Y = N_J^2[z_J']$ and $G_J = G - Y$. Then, $|Y| \geq 5$, $z_J \in V(J - Y)$, $|N_{J-Y}[z_J]| \geq 3$, $G_J$ has a component $G_v^*$ such that $N[v] \cup \{y_I\} \cup N_{J-Y}[z_J] \subseteq V(G_v^*)$, and ${\rm C}(G_J) \subseteq \{G_v^*, (\{9\}, \emptyset)\}$. We have $|V(G_v^*)| \geq 8$, $d_{G_v^*}(v) = d(v) = 3$, $3 = \Delta(G) \geq d_{G_v^*}(z_J) \geq |\{y_I\} \cup N_{J-Y}(z_J)| \geq 3$ (recall that $|N_{J-Y}[z_J]| \geq 3$), and hence $N_{G_v^*}(z_J) = \{y_I\} \cup \{z_J^1, z_J^2\}$ for some distinct $z_J^1$ and $z_J^2$ in $V(J-Y)$. Note that if $i \in [4]$ and $u \in V(G_9^i)$ with $d_{G_9^i}(u) = 3$, then $\Delta(G_9^i - N_{G_9^i}[u]) = 2$. Thus, since $\Delta(G_v^* - N_{G_v^*}[z_J]) = d(v) = 3$, $G_v^*$ is not a $(\mathcal{G}_8 \cup \mathcal{G}_9)$-graph. By Lemma~\ref{lemma}, Lemma~\ref{lemmacomp} and the induction hypothesis, 
$$\iota_2(G,P_3) \leq 1 + \iota_2(G_J,P_3) = 1 + \iota_2(G_v^*,P_3) \leq \frac{|Y|}{5} + \frac{|V(G_v^*)|}{5} \leq \frac{n}{5}.$$

Now suppose $\mathcal{I}_8 \cup \mathcal{I}_9 = \emptyset$. Then, $\mathcal{I}' = \mathcal{I}_3 \cup \mathcal{I}_4$. Suppose $\mathcal{I}_3 \neq \emptyset$. Let $J \in \mathcal{I}_3$. Let $Y_J = N[v] \cup \{y_I, z_J\}$. Then, $Y_J \subseteq N^2[x_I]$ and $|V(J-z_J)| = 2$. By Lemma~\ref{lemma}, Lemma~\ref{lemmacomp} and the induction hypothesis, 
$$\iota_2(G, P_3) \leq 1 + \iota_2(G - Y_J, P_3) \leq 1 + \sum_{H \in \mathcal{H}^* \setminus \{I\}}\frac{|V(H)|}{5} + \sum_{H \in \mathcal{I} \setminus \{J\}}\frac{|V(H)|+2}{5} < \frac{n}{5}$$
as $1 = (|Y_J| + |V(J - z_J)| - 3)/5$ and $|\mathcal{I} \setminus \{J\}| \leq 1$ (recall that $|\mathcal{I}| \leq 2$). Now suppose $\mathcal{I}_3 = \emptyset$. Then, 
$$\mathcal{I}' = \mathcal{I}_4.$$ 

Suppose $|\mathcal{I}_4| = 2$. Then, $\mathcal{I} = \mathcal{I}_4$. Let $Z = N[v] \cup \bigcup_{H \in \mathcal{I}} (N^2[y_I] \cap V(H))$.
Then, $Z \subseteq N^2[\{v, y_I\}]$, $\mathcal{H}^* \setminus \{I\} \subset {\rm C}(G-Z)$ and $|V(H)| \leq 2$ for each $H \in {\rm C}(G-Z) \setminus (\mathcal{H}^* \setminus \{I\})$. By Lemmas~\ref{lemma} and \ref{lemmacomp}, 
$$\iota_2(G, P_3) \leq |\{v, y_I\}| + \iota_2(G-Z, P_3) = 2 + \sum_{H \in \mathcal{H}^* \setminus \{I\}} |D_H| \leq 2 + \sum_{H \in \mathcal{H}^* \setminus \{I\}}\frac{|V(H)|}{5} < \frac{n}{5}$$
as $|N[v] \cup \{y_I\} \cup \bigcup_{H \in \mathcal{I}_4} V(H)| > 10$.

Now suppose $|\mathcal{I}_4| \neq 2$. Recall that $\emptyset \neq \mathcal{I}' = \mathcal{I}_4$ and $|\mathcal{I}| \leq 2$. Thus, $|\mathcal{I}_4| = 1$. Let $J$ be the member of $\mathcal{I}_4$. Let $Z = \{y_I\} \cup V(J)$. Let $G_Z = G-Z$. Then, $G_Z$ has a component $G^*_v$ such that $N[v] \cup \bigcup_{H \in \mathcal{H}^* \setminus \{I\}} V(H) \subseteq V(G^*_v)$, and ${\rm C}(G_Z)$ is $\{G^*_v\}$ or $\{G^*_v, L\}$ with $\{L\} = \mathcal{I}^*$. Since $G[Z]$ is a connected $5$-vertex graph, Lemma~\ref{5vertex} yields $Z \subseteq N^2[u]$ for some $u \in Z$. 
If $G_v^*$ is not a $\mathcal{G}$-graph, then by Lemma~\ref{lemma}, Lemma~\ref{lemmacomp} and the induction hypothesis, 
$$\iota_2(G, P_3) \leq 1 + \iota_2(G_Z, P_3) = \frac{|Z|}{5} + \sum_{H \in {\rm C}(G_Z)} \iota_2(H, P_3) \leq \frac{|Z|}{5} + \sum_{H \in {\rm C}(G_Z)} \frac{|V(H)|}{5} = \frac{n}{5}.$$
Suppose that $G_v^*$ is a $\mathcal{G}$-graph. Since $d_{G_Z}(v) = 3$, $G_v^*$ is a $(\mathcal{G}_4 \cup \{G_9^1, G_9^2, G_9^3, G_9^4\})$-graph. If ${\rm C}(G_Z) = \{G^*_v, L\}$ with $\{L\} = \mathcal{I}^*$ and $|V(L)| = 5t + r$ for some $t, r \in \{0\} \cup \mathbb{N}$ such that $1 \leq r \leq 4$, then by Lemma~\ref{lemma}, Lemma~\ref{lemmacomp} and the induction hypothesis, 
\begin{align} \iota_2(G, P_3) &\leq 1 + \iota_2(G_Z, P_3) = \frac{|Z|}{5} + \iota_2(G_v^*, P_3) + \iota_2(L, P_3) \nonumber \\
&\leq \frac{|Z|}{5} + \frac{|V(G_v^*)|+1}{5} + \left \lfloor \frac{|V(L)|}{5} \right \rfloor \leq \frac{n}{5} \nonumber
\end{align}
as $\floor{|V(L)|/5} = (|V(L)|-r)/5$. Suppose that if ${\rm C}(G_Z) = \{G^*_v, L\}$ with $\{L\} = \mathcal{I}^*$, then $|V(L)| = 5t$ for some $t \in \mathbb{N}$. 

Let $Y_Z = N^2[v] \cap Z$. Then, $y_I \in Y_Z$. Let $z_J^* \in N_J(z_J)$ and $X_1 = \{v, x_I, y_I, z_J, z_J^*\}$. Then, $X_1 \subseteq N^2[y_I]$. 

For any graph $H$, let ${\rm C}(H,P_3)$ denote the set of components of $H$ that contain some $3$-path. Thus, $\iota_2(F,P_3) = 0$ for any $F \in {\rm C}(H) \setminus {\rm C}(H,P_3)$.
\\
\\
\textbf{Case~2.1:} \emph{$G_v^* \simeq G_9^j$ for some $j \in [4]$}. By Observation~\ref{8cycleobs}, $G_v^*$ contains an $8$-cycle $F = (\{v, x_1, \dots, x_7\}, \{vx_1, x_1x_2, \dots, x_6x_7, x_7v\})$ and has a vertex $x_8$ such that $\{x_8\} = V(G_v^*) \setminus V(F)$ and $vx_8 \in E(G_v^*) \setminus E(F) \subseteq \{vx_8, x_1x_8, x_2x_8\}$. Let $X_2 = \{x_2, \dots, x_6\}$. Then, $X_2 \subseteq N^2[x_4]$. 

Suppose $Y_Z = \{y_I\}$. Suppose ${\rm C}(G_Z) = \{G^*_v\}$. Then, ${\rm C}(G-(X_1 \cup X_2),P_3) = \emptyset$. Since $X_1 \cup X_2 \subseteq N^2[\{y_I, x_4\}]$, Lemmas~\ref{lemma} and \ref{lemmacomp} yield $\iota_2(G,P_3) \leq 2 = (|X_1| + |X_2|)/5 < n/5$. Now suppose ${\rm C}(G_Z) = \{G^*_v, L\}$ with $\{L\} = \mathcal{I}^*$. Then, ${\rm C}(G-(X_1 \cup X_2),P_3) = \{L^+\}$ with $V(L) \subseteq V(L^+) \subseteq V(L) \cup (N(v) \setminus \{x_I\})$. Since $|V(L)| = 5t$, $5t \leq |V(L^+)| \leq 5t+2$. Thus, $L^+$ is not a $\mathcal{G}$-graph. By Lemma~\ref{lemma}, Lemma~\ref{lemmacomp} and the induction hypothesis, $\iota_2(G,P_3) \leq 2 + t = (|X_1| + |X_2| + |V(L)|)/5 < n/5$. 

Next, suppose $|Y_Z| \geq 3$. Let $Z_v = N[v] \cup Y_Z \cup X_2$. Then, $Z_v \subseteq N^2[\{v,x_4\}]$ and ${\rm C}(G-Z_v,P_3) = \mathcal{I}^*$. By Lemmas~\ref{lemma} and \ref{lemmacomp}, 
$$\iota_2(G,P_3) \leq 2 + \iota_2(G-Z_v,P_3) =  2 + \sum_{H \in \mathcal{I}^*} \iota_2(H,P_3) < \frac{|Z_v|}{5} + \sum_{H \in \mathcal{I}^*} \frac{|V(H)|}{5} < \frac{n}{5}.$$

Now suppose $|Y_Z| = 2$. Then, $x_I'y_I' \in E(H)$ for some $x_I' \in N(v)$ and $y_I' \in V(J)$, and $Y_Z = \{y_I, y_I'\}$. If $x_I' = x_I$, then $N[v] \cup Y_Z \subseteq N^2[x_I]$, $|N^2[x_I] \cap Z| \geq 3$, and hence $\iota_2(G,P_3) < n/5$ by the argument for the case $|Y_Z| \geq 3$ (with Lemma~\ref{lemma} applied to $\{x_I, x_4\}$ instead of $\{v, x_4\}$). Suppose $x_I \neq x_I'$. Since $J$ is connected, there exist $z_J^1, z_J^2 \in V(J) \setminus \{y_I'\}$ such that $z_J^1 \neq z_J^2$, $z_J^1 \in N_J(y_I')$ and $z_J^2 \in N_J^2[y_I']$. Let $z_J^3$ be the only member of $V(J) \setminus \{y_I', z_J^1, z_J^2\}$.

Suppose first that ${\rm C}(G_Z) = \{G^*_v, L\}$ with $\{L\} = \mathcal{I}^*$. Let $Z_v = N[v] \cup X_2 \cup \{y_I', z_J^1\}$ and $X_3 = \{x_I', x_4\}$. Then, $Z_v \subseteq N^2[X_3]$, $y_Iz_L \in E(G-Z_v)$ and ${\rm C}(G-Z_v, P_3) = \{L^+\}$ with $\{y_I\} \cup V(L) \subseteq V(L^+) \subseteq \{y_I, z_J^2, z_J^3\} \cup V(L)$. Since $|V(L)| = 5t$, $5t+1 \leq |V(L^+)| \leq 5t+3$. Suppose that $L^+$ is a $\mathcal{G}$-graph. Then, $t = 1$, $L^+ \simeq C_8$ and $V(L^+) = \{y_I, z_J^2, z_J^3\} \cup V(L)$. Since $J, L \in {\rm C}(I')$, $N(V(J)) \cap V(L) = \emptyset$. Thus, we have $N_{L^+}(z_J^2) \subseteq \{y_I, z_J^3\}$ and $N_{L^+}(z_J^3) \subseteq \{y_I, z_J^2\}$, contradicting $L^+ \simeq C_8$. Thus, $L^+$ is not a $\mathcal{G}$-graph. By Lemma~\ref{lemma}, Lemma~\ref{lemmacomp} and the induction hypothesis, $\iota_2(G,P_3) \leq |X_3| + \lfloor |V(L^+)|/5 \rfloor = 2 + t < n/5$.

Now suppose ${\rm C}(G_Z) = \{G^*_v\}$. Then, $V(G) = V(G_v^*) \cup Z = V(G_v^*) \cup V(I)$ and $n = 14$. We need to show that $\iota_2(G,P_3) \leq 2$. Since $d(v) = 3$, $\{x_1, x_7\} \cap \{x_I, x_I'\} \neq \emptyset$. By symmetry, we may assume that $x_1 \in \{x_I, x_I'\}$. We have $N^2[v] = N[v] \cup \{y_I, y_I', x_2, x_6\}$. If $x_1 = x_I'$, then $N[v] \cup \{x_2, x_3, y_I', z_J^1\} \subseteq N^2[x_1]$, and hence $\{v, x_1\}$ is a $(P_3,2)$-isolating set of $G$ (recall that $I \in {\rm C}(G-N[v])$, so no vertex of $I$ is adjacent to a vertex of $G_v^* - N[v]$). Suppose $x_1 = x_I$. Then, $N[v] \cup \{x_2, x_3, y_I, z_J\} \subseteq N^2[x_1]$. If $z_J \neq y_I'$, then $\{v, x_1\}$ is a $(P_3,2)$-isolating set of $G$. Suppose $z_J = y_I'$. Then, since $y_Iz_J \in E(G)$, $y_I', z_J^1 \in N^2[y_I]$. We also have $x_1, x_2 \in N^2[y_I]$ and $N[v] \cup \{x_5, x_6\} \in N^2[x_7]$, so $\{x_7, y_I\}$ is a $(P_3,2)$-isolating set of $G$.
\\
\\
\textbf{Case~2.2:} \emph{$G_v^*$ is a $\mathcal{G}_4$-graph}. Then, $V(G_v^*) = N[v]$. By the argument for Case~2.1 (now taking $X_2 = \emptyset$), we may assume that $|Y_Z| = 2$, $x_I'y_I' \in E(H)$ for some $x_I' \in N(v) \setminus \{x_I\}$ and $y_I' \in V(J) \setminus \{y_I\}$, and hence $Y_Z = \{y_I, y_I'\}$. Let $z_J^1$, $z_J^2$ and $z_J^3$ be as in Case~2.1. 

If ${\rm C}(G_Z) = \{G^*_v, L\}$ with $\{L\} = \mathcal{I}^*$, then by the argument for the corresponding case in Case~2.1 (now taking $X_2 = \emptyset$ and $X_3 = \{x_I'\}$), $\iota_2(G,P_3) < n/5$. 

Suppose ${\rm C}(G_Z) = \{G^*_v\}$. Then, $V(G) = V(G_v^*) \cup Z = V(G_v^*) \cup V(I)$ and $n = 9$. The result follows if we show that $G$ has a vertex $z$ such that $\{z\}$ is a $(P_3,2)$-isolating set of $G$. If $y_Iy_I' \in E(G)$, then take $z = y_I$ or $z = y_I'$. Suppose $y_Iy_I' \notin E(G)$. If $x_Iy_I' \in E(G)$, then take $z = x_I$. If $x_I'y_I \in E(G)$, then take $z = x_I'$. Suppose $x_Iy_I', x_I'y_I \notin E(G)$. Recall that $z_J \in N(y_I)$ and $z_J^1 \in N(y_I')$. Let $Z' = \{z_J^1, z_J^2, z_J^3\}$. Since $y_Iy_I' \notin E(G)$, $z_J \in Z'$. If $|N(y_I) \cap Z'| \geq 2$, then take $z = x_I$. Suppose $N(y_I) \cap Z' = \{z_J\}$. Similarly, we may assume that $N(y_I') \cap Z' = \{z_J^1\}$. If $z_J = z_J^1$, then take $z = z_J$. Suppose $z_J \neq z_J^1$. Let $z_J'$ be the only member of $Z' \setminus \{z_J, z_J^1\}$. Since $J$ is connected and $N(y_I') \cap Z' = \{z_J^1\}$, $z_J' \in N(z_J) \cup N(z_J^1)$. We may assume that $z_J' \in N(z_J)$.

Let $x_I''$ be the only member of $N(x) \setminus \{x_I, x_I'\}$. Since $Y_Z = \{y_I, y_I'\}$ and $\Delta(G) = 3$, $N(x_I'') \setminus \{v\}$ is a subset of $\{x_I, x_I', y_I, y_I'\}$ of size at most $2$. If $N(x_I'') \setminus \{v\}$ is one of $\{x_I, y_I'\}$, $\{x_I', y_I'\}$ and $\{y_I, y_I'\}\}$, then take $z = x_I''$. Suppose that it is not. Then, 
\begin{equation} \mbox{$N(x_I'') \setminus \{v\}$ is a subset of one of the sets $\{x_I, x_I'\}$, $\{x_I, y_I\}$ and $\{x_I', y_I'\}$.} \label{x_I''}
\end{equation}

Since $Y_Z = \{y_I, y_I'\}$, $G-N^2[v] = G[Z']$. If $G[Z']$ does not contain a $3$-path, then take $z = v$. Suppose that $G[Z']$ contains a $3$-path. Suppose $z_Jz_J^1 \in E(G)$. Then, $V(G - N^2[y_I]) \subseteq \{x_I', x_I'', y_I'\}$. If ${\rm C}(G - N^2[y_I], P_3) = \emptyset$, then take $z = y_I$. If ${\rm C}(G - N^2[y_I], P_3) \neq \emptyset$, then $x_I'' \in N(x_I') \cup N(y_I')$, and take $z = y_I'$. Now suppose $z_Jz_J^1 \notin E(G)$. Then, $E(G[Z']) = \{z_Jz_J', z_J'z_J^1\}$. Let $F$ be the $8$-cycle with $E(F) = \{vx_I, x_Iy_I, y_Iz_J, z_Jz_J', z_J'z_J^1, z_J^1y_I', y_I'x_I', x_I'v\}$, and let $A = \{x_I''w \colon w \in N(x_I'')\}$. We have $E(G) = E(F) \cup A$. By (\ref{x_I''}), we have $G \simeq G_9^i$ for some $i \in [4]$, a contradiction.~\hfill{$\Box$}

\section{Problems}
We conclude this paper by posing some problems and a conjecture.

We may assume that the vertex set of an $n$-vertex graph is $[n]$. Let 
\[\mathbf{G} = \{G \colon G \mbox{ is a connected graph, } V(G) = [n] \mbox{ for some } n \geq 1\}.\] 
Thus, $\mathbf{G}$ is an infinite set. For any integer $k \geq 0$ and any real number $\alpha > 0$, let 
\[\mathbf{G}_2(k,\alpha) = \{G \in \mathbf{G} \colon \iota_2(G,\mathcal{E}_k) \leq \lfloor \alpha |V(G)| \rfloor\},\] 
and let
\[\mathbf{G}_2(k,\alpha)^* = \{G \in \mathbf{G}_2(k,\alpha) \colon \iota_2(G,\mathcal{E}_k) = \lfloor \alpha |V(G)| \rfloor\} \quad \mbox{and} \quad \mathbf{G}_2[k,\alpha] = \mathbf{G} \setminus \mathbf{G}_2(k,\alpha).\]
Thus, $\mathbf{G}_2[k,\alpha] = \{G \in \mathbf{G} \colon \iota_2(G,\mathcal{E}_k) > \lfloor \alpha |V(G)| \rfloor\}$. In view of Theorems~\ref{result k=0}--\ref{result k=3}, we pose the following problem.

\begin{problem} \label{problem1} (a) Given $k \geq 0$, is there a rational number $c_k$ such that $\mathbf{G}_2{[k,c_k]}$ is finite and $\mathbf{G}_2{(k,c_k)}^*$ is infinite? \\
(b) If $c_k$ exists, then determine $c_k$, $\mathbf{G}_2{[k,c_k]}$ and (at least) an infinite subset of $\mathbf{G}_2(k,c_k)^*$.
\end{problem}
By Theorems~\ref{result k=0}--\ref{result k=3}, $c_0 = 1/3$, $c_1 = 1/4$, $c_2 = c_3 = 1/5$, $\mathbf{G}_2{[0,c_0]} = \{K_1, K_2\}$, $\mathbf{G}_2[1,c_1] = \{G \in \mathbf{G} \colon G \mbox{ is a $\{K_2, P_3, C_3, C_7\}$-graph}\}$, $\mathbf{G}_2[2,c_2] = \{G \in \mathbf{G} \colon G \mbox{ is a $\mathcal{G}$-graph}\}$, $\mathbf{G}_2[3,c_3] = \{G \in \mathbf{G} \colon G \mbox{ is a $(\mathcal{G}_4 \cup \{C_3, C_9\})$-graph}\}$, $\{G \in \mathbf{G} \colon G \simeq B_{n,k+2} \mbox{ for some } n \geq k+3\} \subseteq \mathbf{G}_2{(k,c_k)^*}$ for each $k \in \{0, 1, 2\}$ and $\{G \in \mathbf{G} \colon G \simeq B_{n,4}' \mbox{ for some } n \geq 5\} \subseteq \mathbf{G}_2{(3,c_3)^*}$. 

\begin{conj} For $k \geq 0$, $c_k$ exists.
\end{conj}

Our next problem is stronger than Problem~\ref{problem1}(a).

\begin{problem} Given $k$ and $n$, what is the smallest rational number $c_{k,n}$ such that $\iota_2(G, \mathcal{E}_k) \leq c_{k,n} n$ for every connected $n$-vertex graph $G$?
\end{problem}
By Theorem~\ref{result k=0}, $c_{0,1} = 1$, $c_{0,2} = 1/2$ and $c_{0,n} = \lfloor n/3 \rfloor / n$ for any $n \in \mathbb{N} \setminus \{1, 2\}$. Similarly, Theorems~\ref{result k=1}--\ref{result k=3} give us the values of $c_{k,n}$ for $1 \leq k \leq 3$ and $n \geq 1$.
\\
\\
\noindent
\textbf{Acknowledgements.} Peter Borg was supported by University of Malta Research Seed Fund 2025 grant MATRP14-25.

\end{document}